%
%
%
%
%
\documentstyle{amsppt}
\magnification=\magstep1   
\vcorrection{-1cm}
\input pictex
\NoPageNumbers
\parindent=8mm
\font\gross=cmbx10 scaled\magstep1   \font\abs=cmcsc10
\font\rmk=cmr8  \font\itk=cmti8    \font\ttk=cmtt8

 \def \Aut{\operatorname{Aut}}
 
 \def \sub{{}\subseteq{}}
 \def\sq{\plot 0 0  1 0  1 1  0 1  0 0 /}
 
 \def\E#1{{\parindent=1truecm\narrower\narrower\noindent{#1}\smallskip}}  

 \headline{\ifnum\pageno=1\hfill %
    \else\ifodd\pageno \Rechts \else \Links \fi  \fi}
    \def\Links{\rmk \the\pageno\hfil Schmidmeier\hfil}
    \def\Rechts{\rmk 
        \hfil A Construction of Metabelian Groups\hfil\the\pageno}
%
\noindent\phantom{\rmk [s-metabelian1d, September 22, 2004]}
%
%
        \vglue1truecm
\centerline{\gross A Construction of Metabelian Groups}
        \bigskip
\centerline{By}
        \medskip
\centerline{\abs Markus Schmidmeier}

\footnote""{Mathematics 
Subject Classification (2000): 
Primary 20D10
, Secondary 20K27
, 16G60
.}

        \bigskip\medskip

\E{{\it Abstract.}
In 1934, Garrett Birkhoff has shown that the number of isomorphism classes
of finite metabelian groups of order~$p^{22}$ tends to infinity with $p$.
More precisely, for each prime number $p$ 
there is a family $(M_\lambda)_{\lambda=0,\ldots,p-1}$
of indecomposable and pairwise nonisomorphic metabelian $p$-groups of the
given order.  In this manuscript we use recent results
on the classification of possible embeddings of a subgroup in a finite abelian
$p$-group to construct families of indecomposable metabelian groups,
indexed by several parameters, which have upper bounds on the exponents 
of the center and the commutator subgroup.
}

\document
\medskip\bigskip\noindent
A group $G$ is called {\it metabelian} if its commutator
subgroup $C_1(G)$ is contained in the center $C(G)$; it follows that
the subgroup embedding $C_1(G)\sub C(G)$ 
is an isomorphism invariant of the group $G$.
In [1, \S15], Birkhoff observes that each such invariant can be realized,
and he uses his construction of pairwise nonisomorphic
subgroup embeddings to obtain a family, indexed by one parameter
$\lambda=0,\ldots,p-1$, of pairwise nonisomorphic metabelian $p$-groups.
We recall his result in Theorem~2.  

In fact, the classification of all subgroup embeddings is a problem considered
infeasible.  It is shown in [3] that for $n>6$ the category 
$\Cal S(\Bbb Z/p^n)$ of all embeddings of a subgroup 
in a $p^n$-bounded finite abelian group is of (controlled)
wild representation type and hence admits families of indecomposable
and pairwise nonisomorphic objects which depend on several parameters.
We obtain corresponding statements about families of 
metabelian groups in Theorem~3 and in Corollary~4. 

In these examples and in Birkhoff's, the exponent
of the commutator subgroup is~$p^4$.  
We show in Theorem~5 that this exponent can be reduced 
to $p^3$, at the expense of a higher order and a larger exponent of the group.

\medskip
\centerline{\abs Birkhoff's Construction of Metabelian Groups}

\smallskip
Let $B$ be a finite abelian $p^n$-bounded group and
$A$ a subgroup of $B$ where the
embedding is denoted by $\iota\:A\to B$.
Define the finite group $M=G(A\subseteq B)$ as the semidirect product
$$M=(A\oplus B)\times_\psi D$$
where $D=\Bbb Z/p^m$ is the cyclic group of order equal to
the exponent of $A$, and 
$\psi\:D\to\Aut(A\oplus B)$ is the group map given by
$\psi(d)(a,b)=(a,b+d\iota(a))$.
Note that additive notation is used for the (noncommutative) group
operation. Thus, $M$ is the set $A\oplus B\oplus D$
with the group operation given by
$$(a,b,d)+(a',b',d')\;=\;(a+a',b+d\iota(a')+b',d+d').$$

\smallskip An object $(A\sub B)$
in the category $\Cal S(\Bbb Z/p^n)$ consists of a 
$\Bbb Z/p^n$-module $B$ of finite composition length,
together with a submodule $A$ of $B$. Thus we are dealing with 
embeddings $(A\sub B)$ of a group $A$ in a finite abelian $p^n$-bounded group
$B$.  Morphisms from $(A\sub B)$ to
$(A'\sub B')$ are given by 
the group maps $f\:B\to B'$ which satisfy the condition that $f(A)\sub A'$
holds. The category $\Cal S(\Bbb Z/p^n)$ has the Krull-Schmidt property, so 
every object has a 
decomposition as a direct sum of indecomposable objects; this decomposition
is unique up to isomorphism and reordering. Via the construction~$G$,
families of indecomposable and pairwise nonisomorphic 
objects in $\Cal S(\Bbb Z/p^n)$ give rise
to corresponding families of metabelian groups:

\medskip\noindent{\bf Lemma 1.} 
\item{1.} {\it The group $M=G(A\subseteq B)$ has center 
        $C(M)=B$ and commutator
        subgroup $C_1(M)=\iota(A)$. 
        Hence it is a metabelian group.}
\item{2.} {\it For each $n\in\Bbb N$,
        the map 
        $\Cal S(\Bbb Z/p^n)\to \Cal Groups, (A\subseteq B)\mapsto
            G(A\subseteq B)$, preserves indecomposable 
        objects, and preserves and 
        reflects isomorphisms.}

\smallskip\noindent{\it Proof:\/}
The first assertion is an immediate consequence of the group operation.
If $f\:B\to B'$ gives rise to
an isomorphism $(A\sub B)\to (A'\sub B')$, then 
$A$ and $A'$ have the same exponent, say $p^m$, and the diagonal map
$(f|_{A,A'},f,1_{\Bbb Z/p^m})$ yields an isomorphism 
$G(A\sub B)\to G(A'\sub B')$.
Moreover, the construction $G$ has a left inverse given by assigning to a 
metabelian group $M$ the embedding of the commutator subgroup in the center.
Thus if $(A\sub B)$ is an object in $\Cal S(\Bbb Z/p^n)$ and 
$M=G(A\sub B)$ then $(A\sub B)$ and $(C_1(M)\subseteq C(M))$ are isomorphic
objects in $\Cal S(\Bbb Z/p^n)$. Hence the second assertion follows. \qed

\medskip\noindent{\bf Theorem 2 }(Birkhoff){\bf.} {\it
The number of isomorphism classes of metabelian groups of order $p^{22}$
tends to infinity with $p$. }

\smallskip\noindent{\it Proof:\/}
Let $B$ be the finite abelian $p$-group generated by elements 
$x$, $y$, and $z$ of order $p^6$, $p^4$, and $p^2$, respectively. 
For each value of the parameter $\lambda\in\{0,\ldots,p-1\}$,
let $A_\lambda$ be the subgroup of $B$ generated by $u=p^2x+py+z$
and $v_\lambda=p^2y+p\lambda z$, as indicated in the diagram.
$$
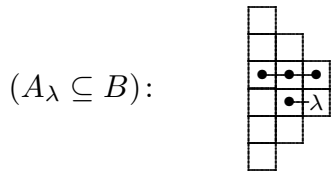
\hbox{\beginpicture 
\setcoordinatesystem units <0.3cm,0.3cm>
\put{$(A_\lambda\sub B)\:$} at -6 2.5
\multiput{\sq} at 0 5  0 4  0 3  0 2  0 1  0 0  1 4  1 3  1 2  1 1  2 3  2 2  /
\put{$\ssize \bullet$} at 0.5 3 
\put{$\ssize \bullet$} at 1.5 3 
\put{$\ssize \bullet$} at 2.5 3 
\put{$\ssize \bullet$} at 1.5 2
\plot 0.5 3  2.5 3 /
\plot 1.5 2  2.2 2 /
\put{$\ssize \lambda$} at 2.5 2
\endpicture}$$
It is shown in [1; Corollary 15.1] that the objects $(A_\lambda\sub B)$
in $\Cal S(\Bbb Z/p^6)$ are indecomposable and pairwise nonisomorphic.
By Lemma 1, the family
$G(A_\lambda\sub B)$, where $0\leq \lambda<p$, consists 
of $p$ indecomposable and pairwise nonisomorphic metabelian
groups; each has order $p^{22}$. 
\qed

\smallskip We observe for later use 
that the metabelian groups constructed this way
all have exponent $p^6$, which is also the exponent of the center, while
the exponent of the commutator subgroup is $p^4$. 
It has been shown in [2] that for $n<6$, there are only finitely many
indecomposable embeddings $(A\sub B)$ in $\Cal S(\Bbb Z/p^n)$, so under 
the above construction, the exponent $p^6$ of the center is minimal.
Conversely, for $n=6$, there are actually infinitely many indecomposable
embeddings in $\Cal S(\Bbb Z/p^n)$, up to isomorphism, and hence 
{\it for each
prime number $p$ there are infinitely many 
indecomposable finite metabelian $p$-groups, up to
isomorphism.}

\medskip
\centerline{\abs Wild Representation Type}

\nopagebreak
\smallskip
In order to construct metabelian groups indexed by several parameters,
we first recall the corresponding situation for subgroup embeddings.
Consider the classical problem of finding a normal form for pairs
of square matrices of the same size with coefficients in some field $k$, 
up to simultaneous conjugation.
This problem is considered infeasible. In particular,
for any given finite number of parameters in $k$, there exist ``families''
of pairwise nonequivalent pairs of square matrices
of the same size, indexed by this parameter set;
here we can even assume that each matrix pair is indecomposable
in the sense
that it is not equivalent to a pair of proper block diagonal matrices 
of the same block type. 
For example, a family indexed by two parameters arises
already in the case of $1\times1$-matrices:  The matrix pairs 
$([\lambda],[\mu])$, where $\lambda,\mu\in k$, 
are all pairwise nonequivalant under simultaneous conjugation.

\smallskip
 Clearly, a pair $(X,Y)$ of 
$m\times m$-matrices defines the structure of a $k\langle X,Y\rangle$-module
on the vector space $V=k^m$ where $k\langle X,Y\rangle$ is the
free algebra in two generators; here the action of the
indeterminates $X$ and $Y$ on $V$ is given by multiplication by the 
matrices with the same names. We denote this module by $(V;X,Y)$. 
It is easy to see that two pairs of
$m\times m$-matrices are equivalent if and only if the two corresponding
$k\langle X,Y\rangle$-modules are isomorphic, so the problem of classifying
all finite dimensional $k\langle X,Y\rangle$-modules, up to isomorphism, is
considered infeasible, too, and hence this classification problem gives
rise to the notion of ``wild type''.

\smallskip An additive category $\Cal A$ is {\it controlled $k$-wild}\/
provided there are full subcategories $\Cal C\sub \Cal B\sub \Cal A$
such that the subquotient $\Cal B/\langle \Cal C\rangle$ of $\Cal A$
is equivalent to the category 
$k\langle X,Y\rangle$-mod.  Here $\langle \Cal C\rangle$ is the categorical
ideal in $\Cal B$ of all morphisms which factor through a (finite) direct
sum of objects in $\Cal C$.  It is shown in [3, Theorem 2] 
that for $n>6$ the category
$\Cal S(\Bbb Z/p^n)$ is controlled $k$-wild where $k$ is the field $\Bbb Z/p$.
In this case, the category $\Cal B$ consists of objects $M$ which are
{\it in between} the objects $I$ and $J$ in the sense that 
$I^\ell\sub M\sub J^\ell$ holds for some $\ell\in\Bbb N$, 
and $\Cal C$ consists of the 
single object $I$.  The two objects $I$ and $J$ are as follows.
$$
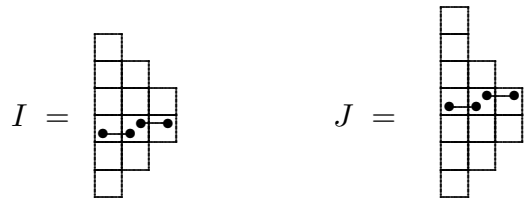
\hbox{\beginpicture 
\setcoordinatesystem units <0.3cm,0.3cm>
\put{$I\;=$} at -2 2.5
\multiput{\sq} at 0 5  0 4  0 3  0 2  0 1  0 0  1 4  1 3  1 2  1 1  2 3  2 2  /
\put{$\ssize \bullet$} at 0.3 1.8 
\put{$\ssize \bullet$} at 1.3 1.8 
\put{$\ssize \bullet$} at 1.7 2.2 
\put{$\ssize \bullet$} at 2.7 2.2 
\plot 0.3 1.8  1.3 1.8 /
\plot 1.7 2.2  2.7 2.2 /
\endpicture}\qquad\quad
\hbox{\beginpicture 
\setcoordinatesystem units <0.3cm,0.3cm>
\put{$\qquad J\;=$} at -4 2.5
\multiput{\sq} at 0 6 0 5  0 4  0 3  0 2  0 1  0 0  1 4  1 3  1 2  1 1  2 3  2 2  /
\put{$\ssize \bullet$} at 0.3 2.8 
\put{$\ssize \bullet$} at 1.3 2.8 
\put{$\ssize \bullet$} at 1.7 3.2 
\put{$\ssize \bullet$} at 2.7 3.2 
\plot 0.3 2.8  1.3 2.8 /
\plot 1.7 3.2  2.7 3.2 /
\endpicture}
$$
Thus, if $I=(I_1\sub I_0)$ and $J=(J_1\sub J_0)$ then the abelian group
$J_0$ is generated by three elements $x,y,z$ of order $p^7,p^4,p^2$,
respectively, and $I_0$ is the subgroup generated by $px,y,z$. 
The subgroup $J_1$ is generated by $p^3x-py$ and $py-z$, and $I_1=pJ_1$.
The equivalence $\Cal B/\langle I\rangle\to k\langle X,Y\rangle$-mod 
is established by providing a full and dense 
functor $F\:\Cal B\to k\langle X,Y\rangle$-mod
and a construction which assigns to the $k\langle X,Y\rangle$-module 
$(V;X,Y)$ the object $S_{(V;X,Y)}$.  If $V=k^m$ then $S_{(V;X,Y)}$ 
is also denoted by $S_{X,Y}$ and satisfies
$I^{2m}\sub S_{X,Y}\sub J^{2m}$.  The functor is left inverse to the
construction in the sense that the $k\langle X,Y\rangle$-modules
$F(S_{(V;X,Y)})$ and $(V;X,Y)$ are isomorphic. For more details on the 
construction (which is obtained by taking fibre products in a 
homomorphism category), 
we refer the reader to the discussion following [3,~Theorem~2]
and the presentation of the construction $\Phi$ after [3,~Proposition~1]. 
As a consequence of the above, the functor $F$ 
preserves indecomposable objects and reflects isomorphisms
(and clearly, as an additive functor, $F$ reflects indecomposable
objects and preserves isomorphisms). 

\medskip
\centerline{\abs Metabelian Groups Indexed by Several Parameters}

\nopagebreak\smallskip
We have seen in Lemma 1 that the construction $G$ preserves indecomposablity
and preserves and reflects isomorphisms.  As a consequence we obtain

\medskip\noindent{\bf Theorem 3.} {\it
Let $m$ be a natural number and let $(X,Y)$ and $(X',Y')$ be two pairs 
of $m\times m$-matrices with 
coefficients in the field $\Bbb Z/p$.

\nopagebreak\smallskip\item{{\rm 1.}} 
The metabelian group $G(S_{X,Y})$ is indecomposable if
the matrix pair $(X,Y)$ is not equivalent under 
simultaneous conjugation
to a pair of proper block diagonal matrices of the same block type.

\item{{\rm 2.}} The two groups $G(S_{X,Y})$ and $G(S_{X',Y'})$ are
isomorphic if and only if the two pairs of matrices are equivalent under
simultaneous conjugation.\qed

}

\smallskip We demonstrate this correspondence 
in the case $m=1$ where we obtain a two parameter
family of indecomposable and pairwise nonisomorphic 
metabelian groups. Given a pair of 
$1\times1$-matrices $([\lambda],[\mu])$, where $0\leq \lambda,\mu <p$,
the corresponding subgroup 
embedding $S_{[\lambda],[\mu]}=(A_{\lambda,\mu}\sub B)$ 
in $\Cal S(\Bbb Z/p^7)$ is as follows.
The abelian group $B$ is generated by elements $x_1$, $x_2$, $y_1$, $y_2$, 
$z_1$, and $z_2$ of order $p^7$, $p^6$, $p^4$, $p^4$, $p^2$, and $p^2$, 
respectively, and
the subgroup $A_{\lambda,\mu}$ 
is generated by $u_1=p^3x_1+py_1+py_2+z_2$,
$u_2(\lambda,\mu)=p^2x_2+p\lambda y_1+p\mu_*y_2+\lambda z_1+\mu z_2$, 
$v_1=p^2y_1+pz_1$, and 
$v_2=p^2y_2+pz_2$; these elements have 
order $p^4$, $p^4$, $p^2$, and $p^2$, respectively, and $\mu_*=\mu+1$.
The embedding $(A_{\lambda,\mu}\sub B)$ can be pictured as follows.
$$
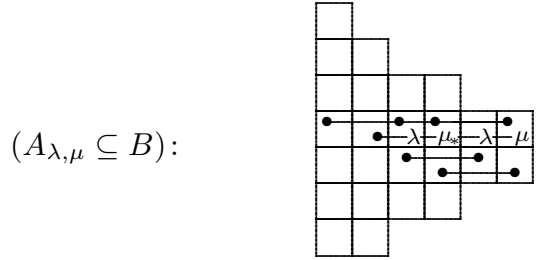
\hbox{\beginpicture 
\setcoordinatesystem units <0.4cm,0.4cm>
\put{$(A_{\lambda,\mu}\sub B)\:$} at -6 2.5
\multiput{\sq} at 0 6  0 5  0 4  0 3  0 2  0 1  0 0  1 5  1 4  1 3  1 2  
        1 1  1 0  2 4  2 3  2 2  2 1  3 4  3 3  3 2  3 1  4 3  4 2  5 3  5 2 /
\multiput{$\ssize \bullet$} at 0.3 3.2  2.3 3.2  3.3 3.2  5.3 3.2  
        1.7 2.8 
        2.5 2.2  4.5 2.2  3.5 1.8  5.5 1.8 /    
\plot 0.3 3.2  5.3 3.2 /
\plot 3.8 3.2  5.2 3.2 /
\plot 1.7 2.8  2.5 2.8 /
\plot 2.9 2.8  3.3 2.8 /
\plot 3.95 2.8 4.5 2.8 /
\plot 4.9 2.8  5.4 2.8 /
\plot 3.5 1.8  5.5 1.8 /
\plot 2.5 2.2  4.5 2.2 /
\multiput{$\ssize \lambda$} at 2.7 2.8  4.7 2.8 /
\put{$\ssize\mu_{\!*}$} at 3.7 2.74 
\put{$\ssize\mu$} at 5.7 2.74   
\endpicture}$$

\smallskip
By applying the construction $G$, we obtain the metabelian group
$G(A_{\lambda,\mu}\sub B)$.  It is generated by 11 elements
$$u_1,\, u_2(\lambda,\mu),\, v_1,\, v_2,\quad x_1,\, x_2,\, y_1,\, y_2,\, z_1, \,z_2,\quad d$$
of order $p^4, p^4, p^2, p^2, \; p^7, p^6, p^4, p^4, p^2, p^2,\; p^4$, 
respectively, and hence $G(A_{\lambda,\mu}\sub B)$ has order $p^{41}$.
The generators commute with each other with the exception of the following
four pairs for which we list the commutators (still using additive notation).
$$\eqalign{ [u_1,d\,]\; &=\; (dp^3)\,x_1 + (dp)\,y_1 + (dp)\,y_2 
                + (d)\,z_2 \cr
        [u_2(\lambda,\mu),d\,]\; &=\; (dp^2)\,x_2 + (d\lambda p)\,y_1 
                + (dp(\mu+1))\,y_2 + (d\lambda)\,z_1 + (d\mu) \,z_2 \cr
        [v_1,d\,]\; &=\; (dp^2)\, y_1 + (dp)\, z_1 \cr
        [v_2,d\,]\; &=\; (dp^2)\, y_2 + (dp)\,z_2 }$$

\medskip We have shown the following result:

\medskip\noindent{\bf Corollary 4.} {\it 
For each prime number $p$, there are at least $p^2$ indecomposable and 
pairwise nonisomorphic metabelian groups of order
$p^{41}$.} \qed

\medskip
\centerline{\abs Commutator Subgroups with Small Exponents}

\nopagebreak\smallskip
In the above example, the exponent of the groups constructed is $p^7$,
while the exponent of the commutator subgroup is $p^4$. 
In the remainder of this note we specify metabelian groups
for which the commutator subgroups all have even smaller exponent~$p^3$;
this is 
at the expense of the order and of the exponent of the group.

\medskip\noindent{\bf Theorem 5.} \nopagebreak
\smallskip\item{1.} {\it For each prime $p$, there are at
least $p$ indecomposable and pairwise nonisomorphic
metabelian groups of order $p^{31}$ and exponent $p^7$
such that the exponent of the commutator subgroup is $p^3$.}
\item{2.} {\it For each prime $p$, there are at
least $\big({p\atop 2}\big)$ indecomposable and pairwise nonisomorphic
metabelian groups of order $p^{60}$ and exponent $p^8$
such that the exponent of the commutator subgroup is $p^3$.}

\smallskip\noindent{\it Proof:\/} 
According to [4], the two families of embeddings 
$(A_\lambda\sub B)_{0\leq\lambda<p}$ and
$(C_{\lambda,\mu}\sub D)_{0\leq \lambda<\mu<p}$
pictured below consist of indecomposable and
pairwise nonisomorphic objects in $\Cal S(\Bbb Z/p^7)$ and in
$\Cal S(\Bbb Z/p^8)$, respectively; clearly the subgroups are bounded by $p^3$.
Now apply Lemma~1. \qed
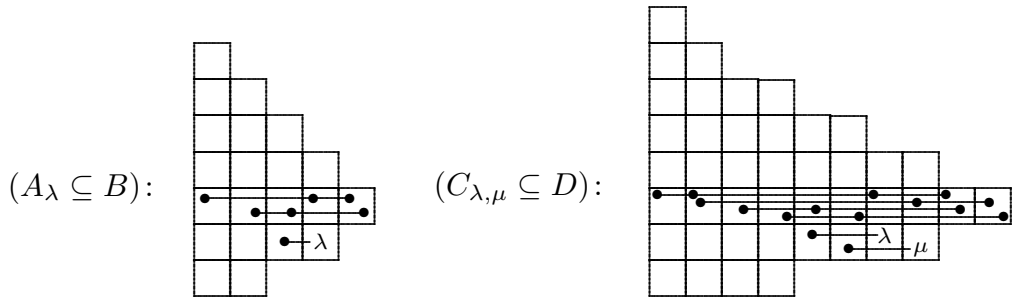
$$\hbox{\beginpicture 
\setcoordinatesystem units <0.4cm,0.4cm>
\put{$(A_{\lambda}\sub B)\:$} at -3 2.5
\multiput{\sq} at 0 6  0 5  0 4  0 3  0 2  0 1  0 0  1 5  1 4  1 3  1 2  
        1 1  1 0  2 4  2 3  2 2  2 1  3 3  3 2  3 1  4 2  /
\multiput{$\ssize \bullet$} at 0.3 2.2  3.3 2.2   4.3 2.2  1.7 1.8  2.7 1.8
        4.7 1.8  2.5 1  /
\plot 0.3 2.2  4.3 2.2 /
\plot 1.7 1.8  4.7 1.8 /
\plot 2.5 1   3.2 1 /
\multiput{$\ssize \lambda$} at 3.5 1 /
\endpicture}
\qquad
\hbox{\beginpicture 
\setcoordinatesystem units <0.4cm,0.4cm>
\put{$(C_{\lambda,\mu}\sub D)\:$} at -3.5 2.5
\multiput{\sq} at 0 7  0 6  0 5  0 4  0 3  0 2  0 1  0 0  
        1 6  1 5  1 4  1 3  1 2  1 1  1 0  2 5  2 4  2 3  2 2  
        2 1  2 0  3 5  3 4  3 3  3 2  
        3 1  3 0  4 4  4 3  4 2  4 1  5 4  5 3  5 2  5 1  
        6 3  6 2  6 1  7 3  7 2  7 1  8 2  9 2 /
\multiput{$\ssize \bullet$} at 0.2 2.3  1.2 2.3  1.4 2.1  6.2 2.3   7.4 2.1  
        8.2 2.3  9.4 2.1  2.6 1.9  3.8 1.7  4.6 1.9  5.8 1.7  
        8.6 1.9  9.8 1.7  4.5 1.2  5.5  .8 /
\plot 0.2 2.3  8.2 2.3 /
\plot 1.4 2.1  9.4 2.1 /
\plot 2.6 1.9  8.6 1.9 /
\plot 3.8 1.7  9.8 1.7 /
\plot 4.5 1.2  6.3 1.2 /
\plot 5.5 .8   7.2 .8 /
\multiput{$\ssize \lambda$} at 6.5 1.2 /
\put{$\ssize \mu$} at 7.5 .8 
\endpicture}$$

        \medskip
\centerline{\abs References}
\baselineskip=9pt \rmk
\smallskip\noindent
\item{[1]} G.~Birkhoff, {\itk Subgroups of Abelian Groups},
Proc.~Lond.~Math.~Soc., II. Ser. 38, 1934, 385--401.

\smallskip\noindent
\item{[2]}  F.\ Richman and  E.\ A.\ Walker, {\itk Subgroups of 
$\ssize p\sp 5$-Bounded Groups,}  Abelian Groups and Modules (Dublin, 1998),  
Trends Math., Birkh\"auser, Basel, 1999, 55--73. 

\smallskip\noindent
\item{[3]} C.~M.~Ringel and M.~Schmidmeier, 
{\itk Submodule Categories of Wild Representation Type},
manu\-script, 2004, 1--13, see {\ttk http://arxiv.org/abs/math/0409417.}

\smallskip\noindent
\item{[4]} M.~Schmidmeier, {\itk Bounded Submodules of Modules,} 
manuscript, 2004, 1--32, preprint available under 
{\ttk http://arxiv.org/abs/math/0408181.}

\medskip\noindent
Markus Schmidmeier, Department of Mathematical Sciences, 
Florida Atlantic University
\par\noindent Boca Raton, Florida 33431-0991, United States of America
\par\noindent {\ttk markus\@math.fau.edu}

\enddocument
\bye